\documentclass[a4paper]{article}

\usepackage{latexsym}
\usepackage{amssymb}
\usepackage{amsmath}
\usepackage{amsthm}
\usepackage{mathrsfs}

\theoremstyle{definition}

\newcommand{\llinf}[2]{\norm{#1}{\L^\infty(#2)}}
\newcommand{\Gl}{\ensuremath{\G_{*,\text{loc}}}}

\renewcommand{\L}{\ensuremath{\mathrm{L}}}
\newcommand{\Con}{\ensuremath{\mathrm{C}}}
\newcommand{\Cinf}{\ensuremath{\mathrm{C}^\infty}}
\newcommand{\D}{\ensuremath{{\cal D}}}
\renewcommand{\S}{\ensuremath{{\cal S}}}

\newcommand{\OM}{\ensuremath{{\cal O}_\mathrm{M}}}

\newcommand{\OC}{\ensuremath{{\cal O}_\mathrm{C}}}
\newcommand{\Zin}{\ensuremath{C_*}}

\newcommand{\mb}[1]{\ensuremath{\mathbb{#1}}}
\newcommand{\N}{\mb{N}}

\newcommand{\R}{\mb{R}}
\newcommand{\C}{\mb{C}}

\newcommand{\G}{\ensuremath{{\cal G}}}

\newcommand{\EM}{\ensuremath{{\cal E}_{\mathrm{M}}}}

\newcommand{\NN}{\ensuremath{{\cal N}}}

\renewcommand{\d}{\ensuremath{\partial}}
\newcommand{\diff}[1]{\frac{d}{d#1}}

\newcommand{\grad}{\ensuremath{\mbox{\rm grad}\,}}

\newfont{\bl}{msbm10 scaled \magstep2}

\newtheorem{thm}{Theorem}
\newtheorem{lemma}[thm]{Lemma}
\newtheorem{prop}[thm]{Proposition}
\newtheorem{defn}[thm]{Definition}

\newtheorem{rem}[thm]{Remark}
\newtheorem{ex}[thm]{Example}

\newcommand{\beq}{\begin{equation}}
\newcommand{\eeq}{\end{equation}}

\newcommand{\col}{\colon}

\newcommand{\FT}[1]{\widehat{#1}}

\newcommand{\F}{\ensuremath{{\cal F}}}

\newcommand{\notmid}{\mid\kern-0.5em\not\kern0.5em}

\newcommand{\norm}[2]{{\| #1 \|}_{#2}}
\newcommand{\lone}[1]{\norm{#1}{L^1}}

\newcommand{\linf}[1]{\norm{#1}{L^\infty}}
\newcommand{\inorm}[2]{\ensuremath{|#1|^*_{#2}}}

\newcommand{\al}{\alpha}
\newcommand{\be}{\beta}
\newcommand{\ga}{\gamma}

\newcommand{\eps}{\varepsilon}

\newcommand{\vphi}{\varphi}

\newcommand{\la}{\lambda}

\newcommand{\Om}{\Omega}
\newcommand{\sig}{\sigma}

\newcommand{\supp}{\mathop{\mathrm{supp}}}
\newcommand{\floor}[1]{\ensuremath{\lfloor #1 \rfloor}}
\renewcommand{\Re}{\ensuremath{\text{Re}}}

\newcommand{\ovl}[1]{\overline{#1}}

\renewcommand{\S}{\mathscr{S}}

\begin{document}

\title{H\"older-Zygmund regularity in algebras of generalized functions}

\author{G\"{u}nther H\"{o}rmann\thanks{Supported by FWF grant P14576-MAT.}\\
Institut f\"{u}r Technische Mathematik, Geometrie und
    Bauinformatik\\ Universit\"{a}t Innsbruck,
    Austria\thanks{Permanent affiliation: Institut f\"{u}r Mathematik,
Universit\"{a}t Wien, Austria} \\
E-Mail: {\tt guenther.hoermann@univie.ac.at}}

\date{July 2003}

\maketitle

\begin{abstract}
We introduce an intrinsic notion of H\"older-Zygmund regularity for 
Colombeau generalized functions. In case of embedded distributions belonging
to some Zygmund-H\"older space this is shown to be consistent. The definition
is motivated by the well-known use of Littlewood-Paley decomposition in
characterizing H\"older-Zygmund regularity for distributions. It is based on a
simple interplay of differentiated convolution-mollification with wavelet
transforms, which directly translates wavelet estimates into properties of
the regularizations. Thus we obtain a scale of new subspaces of the Colombeau
algebra. We investigate their basic properties and indicate first
applications to differential equations whose coefficients are non-smooth but
belong to some H\"older-Zygmund class (distributional or generalized).
In applications problems of this kind occur, for example, in seismology when
Earth's geological properties of fractal nature have to be taken into account
while the initial data typically involve strong singularities. 
\end{abstract}


\section{Introduction}

When studying models of wave propagation in highly irregular media, e.g., in
seismology, (hyperbolic) partial differential equations have to be considered
with coefficients and initial data being generalized functions. The coefficients
represent the medium properties, which may be irregular, e.g., due to folds,
fault zones, or junctions of different geological units as well as caused by long term
physical processes within geological layers. Once the  location
of layer boundaries through jump discontinuities is completed, refined
geological information is reflected in a specific type of regularity patterns of the
material properties within a certain unit. Often self-similar or multi-fractal
behavior can be observed and H\"older continuity, and more generally,
H\"older-Zygmund spaces, were found to be a useful tool for a systematic  
qualitative  analysis (cf.~\cite[Chap.4]{Holschneider:95} and
\cite[Chap.IV]{Triebel:97} for a mathematical justification, and  
\cite{Herrmann:97a,Herrmann:97b,LH:95,SM:93,Wapenaar:98} for seismological applications).

In general, differential equations of the type mentioned above need not make
sense or may fail to have solutions within the theory of distributions. However,
embedding the singular coefficients first into an algebra of generalized
functions, here Colombeau algebras, enables one to carry out a detailed analysis
and yields unique solvability under mild conditions
(cf.~\cite{HdH:01,LO:91,O:89}).

A preliminary study of this procedure in Colombeau
theory was undertaken in \cite{HdH:01c}, where the focus was on 
microlocal properties and the regularization aspects of wavelet transforms. 
The feasibility of recovering Zygmund-H\"older spaces of
positive regularity in one space dimension after the embedding into Colombeau
algebras was proven. In the present paper we extend this result to arbitrary
dimension and regularity scale, although by slightly
changing the definition proposed earlier. We also give first applications to
simple differential equations. In particular, we study a 
(1+1)-dimensional hyperbolic Cauchy problem with typical geophysical conditions on
the coefficients. We show how the regularity of the measured wave depends on
the regularity properties of the medium as well as of the initial value.

The outline of this paper is as follows. After a brief introduction to the
basics of Colombeau theory in Subsection 1.1 we devote Subsection 1.2 to a
review of distributional H\"older-Zygmund spaces and their characterization in
terms of Littlewood-Paley decompositions and wavelet transforms. Section 2
introduces the corresponding  
Colombeau-theoretic notion and discusses basic properties and illustrative
examples. Section 3 presents simple case studies in applications to differential
equations.

\subsection{Colombeau algebras of generalized functions}

We recall the basic facts about the so-called special Colombeau algebras on
$\R^n$. They can be defined on arbitrary open subsets, or even on smooth
manifolds, contain the space of Schwartz distributions, and provide far reaching
consistency with respect to analysis in distribution spaces. For further details
and applications we refer to \cite{Colombeau:85,GKOS:01,O:92}.

The key ingredient of Colombeau algebras is regularization by nets of smooth
functions and the use of asymptotic estimates with respect to the regularization
parameter $\eps$. More precisely, it is based on a quotient construction as
follows: we set (with $I=(0,1]$)
\begin{align*}
        {\cal E} &:= C^\infty(\R^n)^I \\
        \EM &:=\{(u_\varepsilon)_{\varepsilon\in I}\in {\cal E}\mid
                \forall K\subset\subset\R^n, \forall\alpha\in\N_0^n
                \ \exists N\in \N: \\
                &\hspace{3cm}
                \sup_{x\in K}|\partial^\alpha u_\varepsilon
                (x)|=O(\varepsilon^{-N})\mbox{ as }\varepsilon\rightarrow 0\}\\
        \NN &:=\{(u_\varepsilon)_{\varepsilon\in I}\in \EM \mid
                \forall K\subset\subset\R^n, \forall m\in\N:\\
                &\hspace{3.3cm}
                \sup_{x\in K}|u_\varepsilon
                (x)|=O(\varepsilon^{m})\mbox{ as }\varepsilon\rightarrow 0\}.
\end{align*}
$\EM$ is a differential algebras with component-wise operations, $\NN$ is an
ideal in $\EM$, and the {\em special Colombeau algebra} is defined as the
quotient space
\[
        \G :=\EM \,/\,\NN \,.
\]
Since we consider only this type of algebras here we will omit the term
`special' henceforth. A representative of an element $u$ of $\G$ will be denoted
by $(u_\eps)_\eps$, and we will write $u = [(u_\eps)_\eps]$ in this case. Smooth
functions are embedded as a differential subalgebra simply by
$\sigma(f)=[(f)_\eps]$.

To embed nonsmooth distributions we first have to fix a mollifier
$\rho\in\S(\R^n)$ with unit integral satisfying the moment conditions
$\int\rho(x)\,x^\al\,dx = 0$ when $|\al|\geq 1$. Setting
$\rho_\eps(x)=\eps^{-n} \rho(x/\eps)$, compactly supported distributions are
embedded by $\iota_0(w)=(w*\rho_\eps)_\eps+\NN$. Using partitions of unity and
suitable cut-off functions one may explicitly construct an embedding $\iota_\rho
\colon \D' \hookrightarrow \G$ extending $\iota_0$, commuting with partial
derivatives and its restriction to $\Cinf$ agreeing with $\sigma$. Note that
although $\iota_\rho$ depends on the choice of the mollifier $\rho$ this rather
reflects a fundamental property of nonlinear modeling where the interaction of
singular objects depends on the regularization. Additional specifications of the
regularization from a physical model may and should enter the
mathematical theory at this point.

The ring of generalized complex numbers $\widetilde{\C}$ is defined as the
set of moderate nets of numbers ($(r_\eps)_\eps \in \C^I$ with $|r_\eps| =
O(\eps^{-N})$ for some $N$) modulo negligible nets ($|r_\eps| =
O(\eps^{m})$ for each $m$).

\subsection{Review: H\"older-Zygmund regularity of temperate distributions}

This section is a synthesis of related parts from the following sources: in
the basic notation and setup of Zygmund spaces we stay close to
\cite{Hoermander:97}; all wavelet aspects are taken from \cite{Meyer:98};
for further properties of Zygmund classes and related spaces we refer to
\cite{Triebel:I,Triebel:II}.

The result reviewed here is not new and neither are the techniques of its
proof, given in the Appendix. However, we felt the need to unify various
aspects which are crucial to our application later on. The concise summary
of our efforts is the formulation of Theorem \ref{wl_char}.


\paragraph{Continuous Littlewood-Paley decomposition:}


Following \cite[Sect.8.5]{Hoermander:97} we introduce a continuous analog of the
Littlewood-Paley decomposition.

Choose $\vphi\in\D(\R^n)$ real valued and symmetric such that $|\xi| \leq 1$ in
$\supp(\vphi)$ and $\vphi(\xi) = 1$ if $|\xi| \leq 1/2$. Put  $\psi = \diff{t}
\vphi(\xi/t)\mid_{t=1} = - \xi\cdot \grad\vphi(\xi)$ so that the support of
$\psi(./t)$ is contained in the annulus $t/2 \leq |\xi| \leq t$. Observe that we
obtain a continuous partition of unity
 \begin{equation}
    1 = \vphi(\xi) + \int_1^\infty \psi(\xi/t)  \, dt/t.
 \end{equation}

If $f\in\S(\R^n)$ is used as a Fourier multiplier for $u\in\S'(\R^n)$ we will
sometimes write this in pseudodifferential operator notation, i.e.,  $f(D)u =
\F^{-1}(f \FT{u}) = (\F^{-1}f) * u$ where $\F$ and $\FT{\ }$ denote Fourier
transform.

Note that for any $u\in\S'$ and $T \geq 1$ arbitrary we have
 \[
    \vphi(D)u + \int_1^T \psi(D/t) u \, dt/t = \vphi(D/T) u  = T^n
    (\F^{-1}\vphi)(T .) * u
 \]
 which converges to $u$ in $\S'$ when $T\to\infty$.
 This specifies the meaning of the following decomposition formula, which is valid
 in $\S'$,
 \begin{equation}
    u = \vphi(D)u + \int_1^\infty \psi(D/t) u \, dt/t.
 \end{equation}

\paragraph{H\"older-Zygmund spaces:}

The classical H\"older spaces ${\cal C}^s(\R^n)$, for $s > 0$ not integer, as
well as their natural extension to $s\in\N$, the so-called Zygmund classes,
appear in \cite[Section 8.6]{Hoermander:97} in an equivalent realization
given by the spaces $\Zin^s(\R^n)$. These are defined, for any real $s$, in
terms of a continuous Littlewood-Paley decomposition by
 \begin{equation}
    \Zin^s(\R^n) := \{ u\in\S' \mid \inorm{u}{s} := \linf{\vphi(D)u} +
        \sup\limits_{t > 1} \big( t^s \linf{\psi(D/t)u}\big) < \infty \}.
 \end{equation}


Let $m\in\N$. In the context of this paper we call a function
$g\in\S(\R^n)$ a \emph{wavelet of (oscillation) order $m$} if its first $m$ 
moments vanish, that is 
\begin{equation}\label{moments} 
\int x^\al g(x) \, dx = 0 \qquad 0 \leq |\al| \leq m-1
\end{equation}
and it is \emph{weakly radial} (\cite[Chap.\ 1, Equ.\ (5.6)]{Meyer:98}), i.e.,
\begin{equation}\label{w_rad}
    \int_0^\infty |\FT{g}(t\xi)|^2 \frac{dt}{t} = 1 \qquad \forall\xi\not= 0.
\end{equation}
In particular, radial functions can always be normalized so that they
satisfy (\ref{w_rad}).

We introduce the notation $\check{f}(y) = f(-y)$ and $f_\eps(y) =
	\eps^{-n} f(y/\eps)$ for a function $f$ on $\R^n$ (and the bar denoting
complex conjugation).
If $g$ is a wavelet we consider the \emph{wavelet transform} $W_g \colon
\S'(\R^n) \to \Cinf(\R^n\times\R_+)$, mapping $u\in\S'$ into
 \begin{equation}
    W_gu(x,\eps) = u * \ovl{\check{g}_\eps}(x) \qquad
       \forall (x,\eps)\in\R^n\times\R_+.
 \end{equation}
(Note that $\check{}$ denotes reflection, not inverse Fourier transform.)
It is immediate that the image $W_g(\S')$ is
contained in the subspace $\OM(\R^n\times\R_+)$ of smooth functions all of
whose derivatives have polynomial bounds in $x$, $\eps$ and $1/\eps$ (our
notation deviates from \cite{Holschneider:95} where this space is denoted
by $\S'(\R^n\times\R_+)$.) On the space $\OM(\R^n\times\R_+)$ we can define
the \emph{wavelet synthesis operator} $M_g$, mapping $H\in\OM$ into an
element $M_g H \in \S'(\R^n)$, defined by
 \begin{equation} \label{synthesis}
   M_g H = \lim\limits_{r\to 0, R\to\infty}
    \int_r^R H(.,\eps) * g_\eps \frac{d\eps}{\eps}
 \end{equation}
with convergence being understood weakly in $\S'(\R^n)$ (cf.\ \cite[Chapter 1,
Sections 24, 25, and 30]{Holschneider:95}). With the aid of $M_g$ distributions
in $\S'$ can be reconstructed from their wavelet transforms modulo polynomials,
i.e, for each $u\in\S'$ there is a polynomial $p$ on $\R^n$ such that
 \begin{equation} \label{p_inversion}
    u = M_g (W_g u) + p.
 \end{equation}

The crucial observation that motivates the definition of Zygmund regularity
within Colombeau generalized functions is a characterization which is valid for
temperate distributions. As mentioned above this can be found in \cite[Chapter
3]{Meyer:98} in the framework of Bony's two-microlocal spaces. However we repeat
the arguments given there in a `stripped down' version appropriate for the
current context. A detailed proof can be found in the Appendix.

\begin{thm} \label{wl_char}
Let $s$ be a real number and $g\in\S(\R^n)$ have $m$ vanishing moments. Let
$u$ be a temperate distribution on $\R^n$.  
\begin{enumerate}
\item Let $m > s$. If $u \in\Zin^s(\R^n)$ then its wavelet transform satisfies
 \begin{equation} \label{Wg_estimate}
    \sup\limits_{\eps\in(0,1]} \eps^{-s}\, \linf{W_g u (.,\eps)} < \infty.
 \end{equation}
\item Let $m > -s$ and $g$ be weakly radial then (\ref{Wg_estimate}) implies
that there is $u_0\in\Cinf(\R^n)$ such that $u - u_0 \in\Zin^s(\R^n)$.
\end{enumerate}
\end{thm}

\begin{rem} Once more we want to emphasize that the statement of Theorem
\ref{wl_char} is included
in the corresponding, and more general, results presented in Meyer's book
\cite[Chap.~3]{Meyer:98}. The characterization of H\"older-Lipschitz-Zygmund
regularity $0 < s \leq 1$ via the asymptotic behavior of a wavelet-type
transform at small scales has a forerunner in terms of Poisson integrals, e.g.,
in \cite[VII.5]{Zygmund:68} for the one-dimensional case and in
\cite[V.4.2]{Stein:70} on $\R^n$.
\end{rem}

\section{Intrinsic H\"older-Zygmund regularity of \\ Colombeau functions}

\subsection{Basic notions and coherence properties}

We recall that a mollifier is a function $\rho\in\S(\R^n)$ with $\int \rho = 1$.
In addition, we will henceforth assume $\rho$ to be radial.

\paragraph{Mollifiers and wavelets:} We restate the following facts from
\cite[Sect.~3.3]{HdH:01c} \leavevmode
\begin{trivlist}
\item[(i)] Let $\al\in\N_0^n$ with $|\al| \geq 1$. Then the function
$\rho^\al := \ovl{(\d^\al \rho)\check{ }}$ has $m =|\al|$ vanishing moments 
and for any $u\in\S'(\R^n)$
 \begin{equation} \label{mw1}
    \d^\al(u * \rho_\eps)(x) = \eps^{-|\al|} W_{\rho^\al} u (x,\eps)
        \qquad \forall (x,\eps)\in\R^n\times\R_+.
 \end{equation}
In particular, $\ovl{(\Delta^k \rho)\check{\ }}$ is a
wavelet of oscillation order $2k$.
\item[(ii)] If $\int x^\al \rho(x) \, dx = 0$ when $1 \leq |\al| \leq m-1$ then
    $\ovl{\check{\mu}} := - \diff{\eps}(\rho_\eps)\mid_{\eps=1}$ defines
    a wavelet of oscillation order $m$ and for any $u\in\S'(\R^n)$
 \begin{equation} \label{mw2}
    u * \rho_\eps(x) = u * \rho(x) + \int_\eps^1 W_\mu u (x,r) \, \frac{dr}{r}
        \qquad \forall (x,\eps)\in\R^n\times\R_+.
 \end{equation}
\end{trivlist}
In view of Theorem \ref{wl_char} equation (\ref{mw1}) suggests to test for
Zygmund regularity after embedding by looking at the asymptotic properties of
high-order derivatives. The following definition is based on this idea and
refines it in order to ensure mapping properties with respect to
differentiations. Note that it differs from the definition proposed earlier in
\cite{HdH:01c}.
\begin{defn} \label{col_zyg_def} Let $s\in\R$ and $u = [(u_\eps)_\eps] \in\G(\R^n)$.
We say that $u$ is of \emph{(generalized) Zygmund regularity $s$}, denoted
$u\in\G^s_*(\R^n)$, if for $\al\in\N_0^n$ \begin{equation} \label{zyg_def}
    \linf{\d^\al u_\eps} =
    \begin{cases}
        O(1)                 & 0 \leq |\al| < s\\
        O(\log(1/\eps))      & |\al| = s \in\N_0 \\
        O(\eps^{s - |\al|})  & |\al| > s
    \end{cases} \qquad (\eps\to 0).
\end{equation}
\end{defn}

\begin{rem} As a matter of fact, equation (\ref{mw1}) and Theorem \ref{wl_char}
directly suggest to include the third line in (\ref{zyg_def}) of the above
definition. This would be already suitable to characterize the embedded Zygmund
classes (modulo smooth functions) among all embedded temperate distributions as
can be seen from the proof of Theorem \ref{col_zyg_thm} below. However, if we
want the family of spaces $\G_*^s$ ($s\in\R$) to be a scale, in the sense that
$s' \geq s$ implies $\G_*^{s'} \subseteq \G_*^s$, then the testing of decrease
properties must not start at a derivative order which depends on the
(prospective) regularity number. In particular, the case that $s$ is an integer
has to be taken into account, which is done here by the minimum possible, i.e.,
logarithmic, growth rate compatible with embeddings.
\end{rem}

\begin{prop} Let $s$, $s'$, and $r$ be real numbers.
\begin{enumerate}
\item The spaces $\G_*^s$ are nested, i.e., $s' \geq s$ implies
    $\G_*^{s'} \subseteq \G_*^s$.
\item For each $\be\in\N_0^n$ we have a linear map
  $\d^\be \col \G_*^s \to \G_*^{s -|\be|}$.
\item Regularity of products:  $\G_*^r \cdot \G_*^s \subseteq \G_*^p$ where
  $p = r + s$ if  $r,s < 0 $,  $p = \min(r,s)$ if $\max(r,s) > 0$, and
  $p =  \min(r,s)_- $ if $\max(r,s) = 0$. (Here, $c_-$ denotes any number
 $c - \sig$ for $\sig > 0$ arbitrary.)
\end{enumerate}
\end{prop}
\begin{proof}
\emph{Part (i):} If $0 \leq |\al| < s \leq s'$ the assertion is trivial. If
$|\al| = s \leq s'$ then $\linf{\d^\al u}$ is $O(1)$ ($s =s'$) or
$O(\log(1/\eps))$ ($s < s'$), that is $O(\log(1/\eps)$ in any case.

The case $|\al| > s$ leaves us with three subcases for the asymptotic bounds of
$\linf{\d^\al u}$: $s < |\al| < s'$ yields $O(1)$ which is $O(\eps^{s -
|\al|})$; $s < |\al| = s'$ gives $O(\log(1/\eps)$ and hence also $O(\eps^{s -
|\al|}$; finally, in case $|\al| > s' \geq s$ we obtain $O(\eps^{s' - |\al|}$
being again $O(\eps^{s - |\al|}$.

\emph{Part (ii):} We use (\ref{zyg_def}) with $\al$ replaced by $\al + \be$ and
note that $|\al + \be| = |\al| + |\be|$. This gives asymptotic bounds $O(1)$ if
$0 \leq |\al| < s - |\be|$, $O(\log(1/\eps))$ if $|\al| = s - |\be|$, and
$O(\eps^{s - |\be| - |\al|}$ if $|\al| > s - |\be|$.

\emph{Part (iii):} We may assume that $r \leq s$, the opposite case being
completely analogous. Let $u\in\G_*^r$, $v\in\G_*^s$, and $\al\in\N_0^n$. In
estimating $\d^\al (uv)$ we use the Leibniz rule and thus have to find
asymptotic upper bounds for the typical term of the form $\d^\be u_\eps \cdot
\d^{\al -\be} v_\eps$ with $\be\in\N_0^n$ such that $\be \leq \al$. This is done
by combination of the asymptotic growth information about each factor
separately.

If $s < 0$ then the largest growth is due to combinations of the form $O(\eps^{r
- |\be|}) \cdot O(\eps^{s - |\al| + |\be|}) = O(\eps^{r + s - |\al|})$. This
proves the thirst case for the regularity $p$.

If $s = 0$ we only have to check the case $|\al| = 0$ separately. To see this,
note that adding $-\sig$ in the exponents does not decrease the bounds
established above and also captures any occurring logarithmic factors stemming
from $\linf{v_\eps}$. In order $0$ the dominating terms are $O(\eps^r)\cdot
O(\log(1/\eps)) = O(\eps^{r - \sig})$ which proves the second case for $p$.

Finally, let $s > 0$. Assuming $|\al| < s$ implies $|\be| < r$ as well as $|\al
-\be| < r \leq s$ and hence produces only $O(1)$ factors. If $|\al| = s$ then
$|\al - \be| = s$ if and only if $|\be| = 0$ in which case the zero order bound
for $\linf{u_\eps}$ is to be multiplied by $\log(1/\eps)$. Otherwise, i.e., if
$|\be| > 0$, then the factor corresponding to $v$ gives only $O(1)$. It follows
that we obtain the upper bound $O(\log(1/\eps)$ if $r = s > 0$ and $O(\eps^{r -
 s}) = O(\eps^{r - |\al|})$ if $r < s$. If $|\al| > s$ all possible nine
combinations of upper bounds may have to be employed but $O(\eps{r - |\al|})$ is
dominating all of them (since $s > 0$).
\end{proof}

\begin{rem}\leavevmode
\begin{trivlist}
\item[(i)] Compare part (iii) of the Proposition with the distribution
    theoretic result on products in Zygmund spaces
    (cf.~\cite[Prop.8.6.8]{Hoermander:97}): If $u\in\Zin^r$,
    $v\in\Zin^s$ then $u \cdot v$ can be defined (as a weakly
    sequentially continuous bilinear map $\D'\times\D' \to \D'$) if $r
    + s > 0$ and gives an element of Zygmund regularity $\min(r,s)$.
\item[(ii)] We note that the subalgebra $\G^\infty$, defined in
    \cite[Sect.25]{O:92}, reflects a somewhat different concept of
    regularity. First of all, the $\G^\infty$-property is tested on
    compact sets only with $\eps$-asymptotic constant with respect to
    derivative orders but dependent on the compact set.  Furthermore,
    it is easy to give examples of Colombeau functions being very
    regular in one sense but not in the other: if $p$ is a polynomial
    and $\chi$ a smooth cutoff function then the class of $\chi(x)
    p(x/\eps^r)$ is in $\G^\infty$ but it has poor Zygmund regularity
    if $r > 0$; on the other hand, for any $s \in\R$, $\eps^s
    \sin(x/\eps)$ defines a $\G_*^s$-class which is not in
    $\G^\infty$.
\end{trivlist}
\end{rem}

Let $\rho$ be a radial mollifier with all higher moments vanishing. (Hence
$\rho$ can be used to construct wavelets of any oscillation order.) Then we have
the embedding $\iota_\rho\colon \S' \hookrightarrow \G$, $v \mapsto
[(v*\rho_\eps)_\eps]$. We show that under these embeddings the above definition
of the subspaces $\G_*^s \subseteq \G$ is compatible with the distributional
Zygmund classes $\Zin^s$.
\begin{thm} \label{col_zyg_thm} For any $s\in\R$: \leavevmode
\begin{enumerate}
\item  $\iota_\rho(\Zin^s(\R^n)) \subseteq \G^s_*(\R^n)$.
\item If $v\in\S'(\R^n)$ and $\iota_\rho(v)\in\G^s_*(\R^n)$ then there is
    $v_0\in\Cinf(\R^n)$ such that $v - v_0 \in \Zin^s(\R^n)$.
\end{enumerate}
\end{thm}
\begin{proof}

\emph{Part (i):} Let $v\in\Zin^s$ and $\al\in\N_0$. We work through all cases to
be distinguished about the relation of $|\al|$ and $s$.
\begin{description}
\item{$|\al| > s$ and $|\al| > 1$:} Application of (\ref{mw1}) and Theorem
    \ref{wl_char}, (i) (with $m = |\al| > s$) yields
    $ |\d^\al(v * \rho_\eps)(x)| = \eps^{-|\al|} |W_{\rho_\al} v(x,\eps)| =
    \eps^{-|\al|}O(\eps^s)$ ($\eps\to 0$)  uniformly in $x\in\R^n$.
\item{$0 \leq |\al| < s$:} In this case $v\in\Con^{\floor{s}}_b$ and we have
  $ \linf{\d^\al (v * \rho_\eps)} = \linf{(\d^\al v)*\rho_\eps} \leq
    \linf{\d^\al v}  \lone{\rho} = O(1)$,
    where we have used that $\Zin^t \subset \L^\infty$ if $t > 0$
    (\cite[2.3.2/Rem.3]{Triebel:II}).
\item{$|\al| = s \in\N_0$:} Since $\d^\al v\in\Zin^0$ we obtain
  $\linf{W_\mu \d^\al v(.,r)} = O(1)$ in formula (\ref{mw2}) and hence
  \[
    \linf{\d^\al (v * \rho_\eps)} \leq \linf{(\d^\al v) * \rho} + C \int_\eps^1 \frac{dr}{r}
        = O(\log(\frac{1}{\eps})).
  \]
\item{$|\al| = 0 > s$:} Again by (\ref{mw2}) and Theorem \ref{wl_char}, (i),
    noting that $s < 0$, we conclude that
    \[
        |v*\rho_\eps| \leq |v*\rho| + C \int_\eps^1 r^{s-1}\, dr = O(\eps^s).
    \]
\end{description}

\emph{Part (ii):} Choose  $2k > |s|$. Then by (\ref{mw1}) with $\rho^{(2k)} :=
\ovl{(\Delta^k\rho)\check{\ }}$ and applying (\ref{zyg_def}) to $v_\eps = v *
\rho_\eps$ we have
    \[
        \linf{W_{\rho^{(2k)}} v (.,\eps)} = \eps^{2k} \linf{\Delta^k (v_\eps)}
        = \eps^{2k} O(\eps^{s-2k}) = O(\eps^s) \quad (\eps\to 0).
    \]
    The assertion follows from Theorem \ref{wl_char}, (ii) (with $m = 2k >
    |s|$).
\end{proof}

The global $\L^\infty$-bounds used in Definition \ref{col_zyg_def} may be
somewhat too restrictive in certain applications and instead of using a
formulation like `is in $\G_*^s$ modulo a very regular function' we may prefer
to use the following localized version of Zygmund regularity.
\begin{defn} \label{col_zyg_def_loc} Let $\Om \subseteq \R^n$ be open and $s\in\R$.
The Colombeau function $u = [(u_\eps)_\eps] \in\G(\Om)$ is said to be
\emph{locally of generalized Zygmund regularity $s$ in $\Om$}, denoted
$u\in\Gl^s(\Om)$, if for all $K \subset\subset \Om$ and $\al\in\N_0^n$
\begin{equation} \label{zyg_def_loc}
    \llinf{\d^\al u_\eps}{K} =
    \begin{cases}
        O(1)                 & 0 \leq |\al| < s\\
        O(\log(1/\eps))      & |\al| = s \in\N_0 \\
        O(\eps^{s - |\al|})  & |\al| > s
    \end{cases} \qquad (\eps\to 0).
\end{equation}
\end{defn}

\subsection{Examples of regularity under composition}

Let $\Gl^\infty(\Om) = \cap_{s\in\R} \Gl^s(\Om)$ denote the set of functions of
arbitrarily high generalized local Zygmund regularity. In contrast to it we say
that $u$ has no Zygmund regularity, or regularity  $-\infty$, if it is not
contained in $\cup_{s\in\R} \Gl^s(\Om)$.

In the following we will consider the set $\OC(\R^n)$ of smooth functions all of
whose derivatives are of the same polynomial growth, i.e., $u\in\Cinf$ and there
is $M\in\R$ such that for all $\al\in\N_0^n$ we have $|u(x)| = O(|x|^M)$ as $|x|
\to \infty$; in this case, we will say that $u$ is of growth order $M$.

We determine the Zygmund regularity of a simple class of Colombeau functions
obtained by scaling the arguments of smooth functions. These are not obtained by
embedding of distributions and it is a special case of composing a smooth
function with a generalized function. However, nontrivial regularity assertions
about more general cases remain open at this stage. (A very useful result in
$\Zin^s$ spaces concerning composition with smooth functions is
\cite[Prop.8.6.12]{Hoermander:97}.)
\begin{prop}
Let $r$ be a real number.
\begin{enumerate}
\item
Let $f\in\OC(\R)$ of growth order $M\in\R$. Then $u_\eps(x) = f(x/\eps^r)$
defines a Colombeau function $u\in\G(\R)$ which is (at least) of local Zygmund
regularity $s$ if $r < 1$. We have $s = -rM$ if $0 < r \leq 1$ and $M > 0$, $s
= 1-r$ in case $0 < r \leq 1$ and $M \leq 0$, and may put $s = \infty$ when $r
\leq 0$. In general, we have no Zygmund regularity if $r > 1$.
\item Let $p$ be a polynomial of degree $m \not= 0$.
Then $u_\eps(x) = p(x/\eps^r)$ defines a Colombeau function of local Zygmund
regularity $\infty$ if $r \leq 0$. If $r > 0$ we have
\[
    u \in \Gl^s(\R) \Leftrightarrow s \leq - rm.
\]
\end{enumerate}
\end{prop}
\begin{proof} \ \\
\emph{Part (i):} If $r > 1$ we consider the (one dimensional) example $u_\eps(x)
= \sin(x/\eps^r)$. The derivative of order $2k$, evaluated at $x = \pi \eps^r /
2$, gives $\pm \eps^{-2kr}$. But this can never be dominated by $\eps^{s - 2k}$
for all $k\in\N$ and $s$ fixed. Thus $u$ has no Zygmund regularity.

The other extreme case is $r \leq 0$ which always leads to $\eps$-independent
bounds over compact sets in each derivative. Thus we have regularity of
arbitrary order.

We are left with the case $0 < r \leq 1$. Let $\al\in\N_0^n$ then
\[
    |\d^\al u_\eps(x)| = \eps^{- r|\al|} |(\d^\al f)(\frac{x}{\eps^r})|
        \leq C_\al \, \eps^{- r|\al|} (1 + \frac{|x|}{\eps^r})^M.
\]
Let $x$ stay in a fixed compact set. If $M > 0$ the right-hand side is bounded
by some constant times $\eps^{- r (|\al| + M)} = O(\eps^{- r M - |\al|})$.
Finally, if $M \leq 0$ all we can say (in general) is that the right-hand side
is $O(\eps^{- r |\al|}) = O(\eps^{(1 - r)|\al| - |\al|})$ which is $O(1)$ if
$|\al| = 0$ and $O(\eps^{1 - r - |\al|})$ otherwise.

\emph{Part (ii):} The case $r \leq 0$ is obvious since all derivatives have
upper bounds independent of $\eps$ then. So we assume $r > 0$ and note that $p$
is not the zero polynomial since it has degree $ m \geq 1$.

Let $\al\in\N_0^n$ and assume $0 \leq |\al| \leq m$, all higher derivatives
vanish. We have $|\d^\al u_\eps(x)| = \eps^{- r |\al|} |(\d^\al p) (x/\eps^r)|$
which is $\eps^{- r |\al|} O(\eps^{- r (m - |\al|)}) = O(\eps^{-r m})$ if $x$
varies in a compact set. Furthermore, since $\d^\al p$ is a polynomial (nonzero
for some $\al$ of each occurring order) the estimates cannot be improved.

Assume that $u\in\Gl^s$. Since $-rm$ is strictly negative $\eps^{-rm}$ is never
dominated by a constant or logarithmic growth. Hence we have the conditions $s -
k \leq -rm$ when $0 \leq k \leq m$. Setting $k=0$ yields $s \leq -rm$.

On the other hand, $s \leq -rm$ is sufficient to establish the corresponding
Zygmund regularity by the above estimates.
\end{proof}

We end this section with two examples falling into the range of the above
proposition and that further illustrate the different behavior of the notions of
Zygmund- and $\G^\infty$-regularity, in particular, with respect to stability
under smooth compositions.
\begin{ex} \leavevmode
\begin{trivlist}
\item[(i)] We have $v=[(x/\eps)_\eps]\in\G^\infty$ and $v\in\G_*^s
    \Leftrightarrow s \leq -1$ (put $m = r = 1$ in the proposition
    above). Consider the composition $u = \sin \circ\, v$ then $u
    \not\in \G^\infty$ but $u\in\G_*^0$ ($M = 0$, $r = 1$).
\item[(ii)] Similarly, $v = [(1+ x^2/\eps)_\eps]\in\G^\infty$ and
    $v\in\G_*^s \Leftrightarrow s \leq -1$ (use the proposition with
    $m = 2$, $r = 1/2$).  Since $v_\eps \geq 1$ for all $\eps > 0$ we
    may form $u = 1/v \in\G$.  We observe that $u \not\in\G^\infty$:
    at $x = 0$, the values of the derivatives can be read off the
    coefficients in the power series expansion $\sum_k x^{2k}
    (-1)^k/\eps^k$, valid in the interval
    $(-\sqrt{\eps},\sqrt{\eps})$.  From the proposition, with $M =
    -2$, $r = 1/2$, we deduce that $u\in\G_*^{1/2}$.
\end{trivlist}
\end{ex}

\section{Application to linear differential equations with nonsmooth
    coefficients}

\subsection{Solutions with classical H\"older continuity}

We start with the simplest possible case of a differential equation and
mention the well-known elliptic case only briefly. Finally, we sketch how a gain
of regularity can be observed in the hyperbolic case too.

\paragraph{Primitive functions in one dimension:} Let $s$ be any real number and
$u\in\Zin^s(\R)$. If $v\in\D'(\R)$ is a primitive distribution of u, i.e., $v' =
u$, then there is $f\in \Cinf(\R)$ such that
 \begin{equation}\label{prim}
    v - f \in \Zin^{s+1}(\R).
 \end{equation}
To see this we can employ an explicit parametrix of $\diff{x}$, given as
pseudodifferential operator with symbol $h(\xi) = \chi(\xi)/i\xi$ where
$\chi\in\Cinf(\R)$ vanishes near $\xi = 0$ and $\chi(\xi) = 1$ when $|\xi| \geq
1$. (Note that $u\in\S'(\R)$ and $\F((h(D)u)' - u) = (\chi - 1) \FT{u}$ has
compact support; hence $(h(D)u)' - u$ is smooth). It follows that $(v - h(D)u)'
= u - (h(D)u)'$ is smooth and so $v - h(D)u$ must be. But $h(D)$ being of order
$-1$ maps $\Zin^s$ into $\Zin^{s + 1}$ (see \cite[Thm.8.6.14]{Hoermander:97})
which proves (\ref{prim}). Alternatively, we could state that $v$ is locally in
$\Zin^{s+1}$ in the sense that $\vphi v$ belongs to this space for any test
function $\vphi\in\D$.

\paragraph{Elliptic partial differential operators:} Consider $P(x,D) u = f$ where
$P$ is an elliptic partial differential operator of order $m$ with coefficients
and right-hand side $f$ in $\Zin^s$, $s > 0$. Then $u\in\Zin^{s+m}$, i.e., we
observe a gain in regularity by the order of the operator. More precise
statements and related results can be found in \cite[Ch.3]{LU:68}, a concise
summary is \cite[Thm.17.1.1']{Hoermander:V3}.

\paragraph{The embryonic hyperbolic case:} As a resemblance of more realistic models from
geophysics we consider the Cauchy problem
 \begin{equation}\label{Zyg_CP}
    \d_t u + a(x) \d_x u  = 0,   u \mid_{t=0} = b
 \end{equation}
where $a\in\Zin^s(\R)$, $0 < s < 1$, and $b\in\Zin^{s+1}(\R)$. In addition, we
make the following strong positivity and boundedness assumption on the
coefficient: there exist constants $c_1$, $c_2$ such that
 \begin{equation} \label{co_ass}
    0 < c_1 \leq a(x) \leq c_2 \qquad \text{for all } x\in\R.
 \end{equation}
This condition is justified, e.g., if $a$ is of the nature of sound speed in a
certain medium or fluid.

The Cauchy problem (\ref{Zyg_CP}) is easily solved by the method of
characteristics. We point out that, by continuity and positivity of the
coefficient $a$, the characteristic ODE has indeed a unique $\Con^1$ solution.
To make this more explicit we define $A(x) = \int_0^x dr / a(r)$.
Note that $A$ is $\Con^1$, strictly monotone, and that $|A(x)| \leq |x| / c_1$.
Then we set
 \begin{equation}\label{u_sol}
    u(x,t) = b(A^{-1}(A(x) - t))
 \end{equation}
which is directly checked to be the $\Con^1$ solution of (\ref{Zyg_CP}). As an
introduction to the subject of the following two sections we investigate its
H\"older-Zygmund regularity in some detail.
\begin{prop} Let $u$ be the solution of (\ref{Zyg_CP}) given by (\ref{u_sol}). Then the
first order derivatives of $u$ are H\"older continuous of order $s$.
\end{prop}
\begin{proof}
Note that $1/a$ is in $\Zin^s$ which can be seen directly or, alternatively, be
deduced from \cite[Prop.8.6.12]{Hoermander:97} since $a$ is bounded away from
zero. We proceed straightforward in two steps.

The function $h(x,t) = A^{-1}(A(x) - t)$ clearly is $\Con^1$. We first show that
its first order derivatives are H\"older continuous with exponent $s$. We have
$\grad h(x,t) = a(h(x,t))\cdot (1/a(x), -1)$, which is bounded, and
\begin{multline*}
    | \grad h(x,t) - \grad h(y,r)| \\
    \leq
        |a(h(x,t))| |(\frac{1}{a(x)}-\frac{1}{a(y)},0)| +
        |a(h(x,t)) - a(h(y,r))| |(\frac{1}{a(y)}, -1)| \\
    \leq C (|x - y|^s + |h(x,t) - h(y,r)|^s) \leq C' (|x - y|^s + |(x-y,t-r)|^s)
\end{multline*}
with generic constants $C$, $C'$ depending on $a$ only. Hence $\grad h$ is
H\"older continuous of order $s$.

The second step is the composition with $b$. We have $\grad u = b'(h) \cdot
\grad h$ and therefore obtain
\begin{multline*}
    | \grad u(x,t) - \grad u(y,r)| \\
    \leq
        |b'(h(x,t))| |\grad h(x,t) - \grad h(y,r)| +
        |b'(h(x,t)) - b'(h(y,r))| |\grad h(y,r)| \\
   \leq C (|x - y|^s + |h(x,t) - h(y,r)|^s) \leq C' (|x - y|^s + |(x-y,t-r)|^s)
\end{multline*}
where we have used the H\"older continuity, as well as the boundedness, of
$\grad h$ and $b'$.
\end{proof}

\subsection{Primitive functions and a linear first order ODE}

The simplest inhomogeneous case is that of primitive functions in one dimension.
Unlike primitive distributions, a Colombeau primitive function need not gain
regularity, as the following examples illustrate.
\begin{ex} The generalized constants  $[(1/\eps^r)_\eps]$, $r > 0$,
do not have Zygmund regularity higher than $-r$ but nevertheless are primitive
functions of $0$. As a consequence, any Colombeau function allows for primitive
functions with Zygmund regularity arbitrarily low. Furthermore, all primitive
functions of $[(1/\eps^r)_\eps]$ are of the form $[(x/\eps^r)_\eps] + c$ where
$c$ is any generalized constant. The latter can never be of Zygmund regularity
higher than $-r$ thereby showing the existence of Colombeau functions possessing
no primitive function of any higher regularity.
\end{ex}

However, saving a minimum of the classical intuition, we can still control the
regularity of primitive functions obtained from embedded distributions via
integration.
\begin{prop} Let $u\in\iota_\rho(\Zin^s(\R))$, $x_0\in\R$ arbitrary, and define a
primitive function $v$ by the representative $v_\eps(x) = \int_{x_0}^x u_\eps(y)
\, dy$. Then $v$ belongs to $\Gl^{s+1}$.
\end{prop}
\begin{proof}
There is $u_0\in\Zin^s$ such that $u = \iota_\rho(u_0)$. By (\ref{prim}) we can
find $w\in\Zin^{s+1}$ of $u_0$ such that we have $\iota_\rho(u_0) =
\iota_\rho(w') + \sig(g)$ for some smooth function $g$. Hence there is
$(n_\eps)_\eps \in \NN$ such that
\[
    v_\eps(x) = \int_{x_0}^x w'*\rho_\eps(y)\, dy + \int_{x_0}^x
      f(y)\, dy + n_\eps(x).
\]
We observe that, in general, any derivative of order $l \geq 1$ has the asserted
asymptotic estimates since $v_\eps^{(l)}(x) = u_\eps^{(l-1)}$, so only the zero
order estimate has to be investigated separately.

Using $\int_{x_0}^x w'*\rho_\eps(y)\, dy = w*\rho_\eps(x) - w*\rho_\eps(x_0)$ we
obtain, for any compact interval $I$ containing $x$, $x_0$, and of length $|I|$,
\[
    |v_\eps(x)| \leq 2 \sup\limits_{y\in I} |w*\rho_\eps(y)| +
        \sup\limits_{y\in I} \big( |I| |f(y)| + |n_\eps(y)| \big).
\]
The second term on the right-hand side is $O(1)$ on compact subsets with respect
to $x$. Finally, since $w\in\Zin^{s+1}$ we deduce from Theorem \ref{col_zyg_thm}
the required growth properties, according to regularity $s+1$, of the complete
expression.
\end{proof}

In the proposition to follow we give a lower bound for the regularity of the
solution to a linear homogeneous ODE with coefficient from a generalized Zygmund
class. We will impose an additional boundedness condition on this coefficient
and recall: $v\in\G$ is said to be \emph{locally bounded} if $\forall K
\subset\subset \R^n$ there is $C$, $\eps_0
> 0$ such that $\sup_{x\in K} |v_\eps(x)| \leq C$ for all $0 < \eps < \eps_0$.
\begin{prop} Assume $s \geq -1$ and let $a\in\Gl^s(\R)$
  such that $\Re(a)$ is locally bounded. Let $b$ be a generalized
  constant which, considered as a generalized function, is of
  generalized Zygmund regularity $t$ ($t\in\R$). Then
the unique solution $u\in\G(\R)$ to the initial value problem \begin{equation}
\label{ODE}
    \diff{x} u(x) = a(x) u(x), \quad  u(0) = b
\end{equation} belongs to $\Gl^r(\R)$ with $r = s + 1$ if $t > 0$ and
$r = t$ if $t < 0$. When $t = 0$ we have $r = 0_-$ if $s = -1$ and $r
    = 0$ if $s > -1$. Here, $0_-$ stands for any negative number,
    arbitrarily close to $0$.
\end{prop}
\begin{proof}
Existence and uniqueness of the solution $u$ follows from \cite{HO:99}. A
representative is given by
\[
    u_\eps(x) = b_\eps e^{\int_0^x a_\eps(y) \, dy}
\]
where $(b_\eps)_\eps$ is a representative of $b$. By our assumption on
$a$ we have on any compact set $K$
\[
    \llinf{u_\eps}{K} = O(|b_\eps|) \quad (\eps \to 0).
\]
To find sharp asymptotic bounds for the derivatives we first investigate their
algebraic structure. The following assertion is easily proved using the ODE
itself and induction on the derivative order $k$. Let $k\in\N$ then
$u_\eps^{(k)}$ is a linear combination of terms of the following form: with
$m\in\N$, $1 \leq m \leq k$, and $\la\in\N_0^m$ such that $|\la| = k - m$ we
have the expression
\begin{equation} \label{der_struc}
    u_\eps \cdot \prod_{j=1}^m a_\eps^{(\la_j)}.
\end{equation} As noted above the first factor, $u_\eps$, is $O(|b_\eps|)$, so we focus on the product of
derivatives of $a_\eps$.

\emph{Claim:} for any $s \geq -1$ we have, with the notation as in
(\ref{der_struc}),
\begin{equation} \label{claim}
    \prod_{j=1}^m \llinf{a_\eps^{(\la_j)}}{K} =
        \begin{cases}
            O(1)                & k < s + 1\\
            O(\log(1/\eps))     & k = s + 1\\
            O(\eps^{s + 1 - k}) & k > s + 1
        \end{cases}
\end{equation} on compact sets with respect to $x$.
\begin{trivlist}
\item[$\bullet$] If $s < 0$ then each $\la_j \geq 0 > s$ and hence we have the
asymptotic bound $O(\eps^{m s - |\la| }) = O(\eps^{m (s + 1) - k}) = O(\eps^{s
+ 1 - k})$.
\item[$\bullet$] If $s = 0$ let $n$ be the number of $j$'s such that $\la_j =
0$. Then we have the asymptotic upper bound involving $(\log(1/\eps))^n
\eps^{-|\la|} = (\log(1/\eps))^n \eps^{m - k}$. When $m \geq 2$ the second
factor is $O(\eps^{1 - k } \eps)$ where $\eps$ can compensate for the
logarithmic terms. Hence we have a bound $O(\eps^{1 - k})$. When $m = 1$ we
obtain $O(\log(1/\eps))$ if $k = 1$ and $O(\eps^{1 - k})$ otherwise (since $n
= 0$ then).
\item[$\bullet$] Finally, we have to consider the case $s > 0$. We have to
further distinguish three subcases for the relation between $k$ and $s + 1$.

Subcase $k < s + 1$: Since $|\la| = k - m \leq k - 1 < s$ we have that each
$\la_j < s$ and hence an upper bound $O(1)$.

Subcase $k = s + 1$: Now $|\la| \leq s$ and for at most one $j$ we have $\la_j =
s$, all others are less than $s$; hence we obtain an estimate of the form
$O(\log(1/\eps)$.

Subcase $k > s + 1$: Denote by $n$ the number of $j$'s such that $\la_j = s$ and
define $N' := \{j \mid \la_j > s \}$, $n' := |N'|$ (the cardinality of $N'$).
Put $\la'_j = 0$ if $j\not\in N'$ and $\la'_j = \la_j$ otherwise. The asymptotic
upper bound in question is now expressible as $O((\log(1/\eps))^n \eps^{n' s -
|\la'|} )$. If $n' = |\la'| = 0$ this clearly is $O(\eps^{s + 1 - k })$, so we
may assume that $n' \geq 1$. Since $k - m = |\la| \geq n s + |\la'|$ we obtain
$\eps^{-|\la'|} \leq \eps^{m + n s - k}$. Inserting this into the expression for
the asymptotic upper bound we arrive at
$O((\log(1/\eps)^n \eps^{n s}) O(\eps^{n' s + m - k})$.
Here, the first factor is $O(1)$ since $s > 0$ and the second factor is
$O(\eps^{s + 1 - k})$, due to $n' \geq 1$ and $m \geq 1$, as claimed.
\end{trivlist}

Now we come back to (\ref{der_struc}) and use the information from (\ref{claim}).
If $t < 0$ the order zero estimate implies $r \leq t$. By combining
(\ref{der_struc}) with (\ref{claim}) we see that regularity $r = t$ can indeed
be established for any $s \geq -1$.

If $t = 0$ the order zero estimate is logarithmic, due to $|b_\eps|$, and forces
$r \leq 0$. If $s + 1 > 0$ it is seen from (\ref{claim}) that $r = 0$ holds. In
case $s = -1$ and $k > 0$ we have to cope with appearing upper bounds of the
form $O(\log(1/\eps) \eps^{-k})$. This requires subtraction of an arbitrary
small, but still positive, number $\sig$ in the exponent to incorporate the
additional logarithmic factor. (Compare with the situation in the general
multiplication result.)

Finally, if $t > 0$ the factor $|b_\eps| = O(1)$ and the regularity $r = s + 1$
is established directly from (\ref{claim}).
\end{proof}

\begin{rem}\leavevmode
\begin{trivlist}
\item[(i)] Note that if $k$ is very large in (\ref{der_struc}) it may
happen that each $\la_j > s$. In this case, a general upper bound will
be of the form $O(\eps^{m s - |\la|}) = O(\eps^{m (s + 1) - k
})$. Since $m$ may also become arbitrarily large this indicates that
the condition $s + 1 \geq 0$ cannot be dropped in general while
expecting Zygmund regularity of the solution. This is illustrated by
the constant coefficient problem with $a_\eps(x) = i/\eps^r$, $b = 1$
and $r > 0$. The solution (representative) is then $\exp(ix/\eps^r)$,
a sort of `standard counter example' in Colombeau regularity theory.
\item[(ii)]  The boundedness condition on the real part of the coefficient
  cannot be dropped.
Indeed, this can be seen from the constant coefficient problem with $a_\eps(x) =
\log(1/\eps)$ and $b=1$. A Colombeau solution representative is given by
$\exp( x \log(1/\eps))$ which is not Zygmund-regular: The $\L^\infty$-norm taken
over a compact set $K$ grows like $\eps^{-m(K)}$ if $m(K)$ denotes the maximum
of $K$.
\end{trivlist}
\end{rem}

\subsection{A linear hyperbolic Cauchy problem}

As we have indicated in the introduction, if we think of modeling seismic wave
propagation we may encounter fractal-like variations in sound speed, for
example. By the very nature of coefficients representing physical observables
like sound speed, density, elasticity tensors, we see that a positivity
condition on the coefficient(s) is not artificial. We state a first regularity
result for a simple model of this type in one space dimension. It fits nicely
with the classical embryonic case discussed in Subsection 4.1.
\begin{thm} \label{PDE_reg} Let $a = [(a_\eps)_\eps] \in\Gl^s(\R)$, $s
  \geq 0$, and assume there
are positive constants $c_1$, $c_2$ such that $c_1 \leq a_\eps(x) \leq c_2$ for
all $x\in\R$ and $\eps\in(0,1)$. Let $t$ be a real number and $b\in\Gl^t(\R)$.
If $u$ is the (unique) solution of the Cauchy problem
\beq
    \d_t u + a(x) \d_x u  = 0, \quad   u(0)  = b
\eeq
then $u\in\Gl^r(\R^2)$ with $r = \min(t,1)_-$  if $s = 0$, and $r =
\min(t,s+1)$ if  $s > 0$. (As above, $\min(t,1)_-$ denotes any number
approximating $\min(t,1)$ from below.)
\end{thm}
\begin{proof}
We have to determine asymptotic upper bounds of all derivatives of $u_\eps(x,t)
= b_\eps(A_\eps^{-1}(A_\eps(x) - t))$ on compact sets.

We first note that the assumptions on $a$ imply that $h_\eps(x,t) =
A_\eps^{-1}(A_\eps(x) - t)$ maps a compact subset $K$ of $\R^2$ into a fixed
compact subset $K'$ of $\R$, independently of $\eps$. Therefore when doing
estimates on $K$ we may essentially ignore the argument $h_\eps(x,t)$ whenever
appearing as inner function in compositions and write instead the supremum over
$K'$. However, the chain rule will bring out derivatives of $h_\eps$ as
additional factors.

Thus the order zero estimate for $u_\eps$ is simply
 \begin{equation} \label{order_0}
    \llinf{u_\eps}{K} \leq \llinf{b_\eps}{K'}.
 \end{equation}

In the following, let $\al\in\N_0^2$ such that $|\al| \geq 1$.

As a preparation we have to investigate the structure of the higher order
derivatives of $u_\eps$. To simplify notation we drop the subscript $\eps$ in
doing this algebra.

\emph{Claim 1:} $\d^\al u$ is a linear combination of terms of the following
form:
 \begin{equation} \label{der_struc2}
    b^{(l)}(h(x,t)) \cdot \prod\limits_{i=1}^m a^{(\la_i)}(h(x,t))
        \cdot \prod\limits_{j=1}^n a^{(\mu_j)}(x) / a^k(x)
 \end{equation}
where $1 \leq l \leq m \leq |\al|$, $0 \leq n \leq |\al|$, $0 \leq k \leq
|\al|$,  $|\la| = m - l$, and $|\mu| = |\al| - m$, with the notation $\la :=
(\la_i)_{i=1}^m$ and $\mu := (\mu_j)_{j=1}^n$.

We prove (\ref{der_struc2}) by induction on $|\al|$. Concerning the inner
derivatives when applying the chain rule we note that, by definition of $A$, we
have $\d_t h(x,t) = - a(h(x,t))$ and $\d_x h(x,t) = a(h(x,t)) / a(x)$.

The base cases correspond to the first order derivatives $\d_t u(x,t) = -
b'(h(x,t))\cdot a(h(x,t))$ and $\d_x u (x,t) = b'(h(x,t)) a(h(x,t)) / a(x)$, both
complying with the structure of (\ref{der_struc2}).

Assume the claim to be proven already for $|\al|$ and let $\be\in \N_0^2$ with
$|\be| = |\al| + 1$. We distinguish the two cases $\be = \al + e_1$ and $\be =
\al + e_2$ ($e_j$ denoting the standard unit vector in direction $j$).

Case $\be = \al + e_1$: By the induction hypothesis, $\d^\be u = \d_x (\d^\al
u)$ is a linear combination of terms
 \[
    \d_x\Big( b^{(l)}(h(x,t)) \cdot \prod\limits_{i=1}^m a^{(\la_i)}(h(x,t))
        \cdot \prod\limits_{j=1}^n a^{(\mu_j)}(x) / a^k(x)  \Big).
 \]
Application of the Leibniz and chain rules yields four types of terms.\\
Type 1 is
 \[
    b^{(l+1)}(h(x,t))\, a(h(x,t))\, a(x)^{-1} \cdot \prod\limits_{i=1}^m a^{(\la_i)}(h(x,t))
        \cdot \prod\limits_{j=1}^n a^{(\mu_j)}(x) / a^k(x)
 \]
which matches the claim with new quantities $l+1$, $m+1$, $k+1$, and $\la_{m+1} := 0$\\
Type 2, for any $1 \leq r \leq m$, is
 \[
    b^{(l)}(h(x,t)) \cdot \prod\limits_{i \not= r, i=1}^m a^{(\la_i)}(h(x,t))
        \cdot a^{(\la_r + 1)}(h(x,t))\cdot \frac{a(h(x,t))}{a(x)^{-1}}
        \cdot \prod\limits_{j=1}^n \frac{a^{(\mu_j)}(x)}{a^k(x)}
 \]
and satisfies (\ref{der_struc2}) with $k+1$, $m+1$, $\la_r + 1$,
and $\la_{m+1} := 0$ instead.\\
Type 3, for any $1 \leq r \leq n$, is
 \[
    b^{(l)}(h(x,t)) \cdot \prod\limits_{i=1}^m a^{(\la_i)}(h(x,t))
        \cdot \prod\limits_{j \not= r, j=1}^n a^{(\mu_j)}(x)
        \cdot a^{(\mu_r + 1)}(x) / a^k(x)
 \]
where we may use the new component $\mu_r + 1$ in (\ref{der_struc2}).\\
Type 4 is
 \[
    b^{(l)}(h(x,t)) \cdot \prod\limits_{i=1}^m a^{(\la_i)}(h(x,t))
        \cdot \prod\limits_{j=1}^n a^{(\mu_j)}(x) \cdot (-k a'(x)) / a^{k+1}(x)
 \]
and matches the claim with new quantities $k+1$, $n+1$, and $\mu_{n+1} := 1$.

Case $\be = \al + e_2$: By the induction hypothesis, $\d^\be u = \d_t (\d^\al
u)$ is a linear combination of terms
 \[
    \d_t \Big( b^{(l)}(h(x,t)) \cdot \prod\limits_{i=1}^m a^{(\la_i)}(h(x,t))
        \cdot \prod\limits_{j=1}^n a^{(\mu_j)}(x) / a^k(x)  \Big).
 \]
Application of the Leibniz and chain rules yields two types of terms.\\
Type 1 is
 \[
   - b^{(l+1)}(h(x,t))\, a(h(x,t)) \cdot \prod\limits_{i=1}^m a^{(\la_i)}(h(x,t))
        \cdot \prod\limits_{j=1}^n a^{(\mu_j)}(x) / a^k(x)
 \]
which matches the claim with new quantities $l+1$, $m+1$, and $\la_{m+1} := 0$\\
Type 2 is
 \[
    - b^{(l)}(h(x,t)) \cdot \prod\limits_{i \not= r, i=1}^m a^{(\la_i)}(h(x,t))
        \cdot a^{(\la_r + 1)}(h(x,t))\, a(h(x,t))
        \cdot \prod\limits_{j=1}^n a^{(\mu_j)}(x) / a^k(x)
 \]
and satisfies (\ref{der_struc2}) with new values $m+1$, $\la_r + 1$, and
$\la_{m+1} := 0$.\\
The claim is proved.

According to claim 1 and the remark at the beginning of this proof we deduce
that on any compact set we have
 \begin{equation} \label{order_large}
    \llinf{\d^\al u_\eps}{K} = O( \llinf{b_\eps^{(l)}}{K'} \cdot \prod\limits_{i=1}^m
    \llinf{a_\eps^{(\la_i)}}{K'}
        \cdot \prod\limits_{j=1}^n \llinf{a_\eps^{(\mu_j)}}{K} ).
 \end{equation}
To evaluate this carefully we first focus on all the factors having bounds
depending on $a$ or its derivatives.

With the notation of (\ref{der_struc2}) define the sets $L_0 = \{ i \mid \la_i =
s \}$, $L_1 = \{ i \mid \la_i > s \}$, $M_0 = \{ j \mid \mu_j = s \}$, and $M_1
= \{ j \mid \la_j > s \}$. Let $l_0 = |L_0|$ and define similarly $l_1$, $m_0$,
$m_1$ as the respective cardinalities.

\emph{Claim 2:}  On compact sets we can give asymptotic upper bounds of the
following form
 \begin{multline} \label{prod_bound}
    \prod\limits_{i=1}^m \llinf{a_\eps^{(\la_i)}}{K'}
        \cdot \prod\limits_{j=1}^n \llinf{a_\eps^{(\mu_j)}}{K}  =\\
        =
    \begin{cases}
        O((\log(1/\eps))^{l_0 + m_0} \eps^{(l_0 + m_0) s} \eps ^{(l_1 + m_1) s + l - |\al|})
            & l_1 + m_1 > 0\\
        O((\log(1/\eps))^{l_0 + m_0}) & l_1 + m_1 = 0.
    \end{cases}
 \end{multline}
Using the notation introduced above the proof is easy. We observe that each
$i\in L_0$ and $j \in M_0$ contributes a factor $\log(1/\eps)$, whereas each $i
\in L_1$, resp.\ $j \in M_1$, gives rise to a factor $\eps^{s - \la_i}$, resp.\
$\eps^{s - \mu_j}$. We define the tuples $\la'$, resp.\ $\mu'$, by setting all
components in $\la$, resp.\ $\mu$, which are less than $s$ to $0$. Then we
obtain a total bound
 $O((\log(1/\eps))^{l_0 + m_0} \eps^{(l_1 + m_1) s - |\la'| - |\mu'|})$.
If $l_1 + m_1 = 0$ then also $|\la'| + |\mu'| = 0$ which proves the second case
in (\ref{prod_bound}). If $l_1 + m_1 \geq 1$ we note that $m - l = |\la| \geq
|\la'| + l_0 s$ and $|\al| - m = |\mu| \geq |\mu'| + m_0 s$. This implies
$-|\la'| - |\mu'| \geq (l_0 + m_0) s + l - |\al|$ and hence $\eps^{-|\la'| -
|\mu'|} \leq \eps^{(l_0 + m_0) s + l - |\al|}$. Inserting this into the above
total bound matches the first case in (\ref{prod_bound}) and proves claim 2.

We are now in a position to estimate the regularity $r$ of $u$ using
(\ref{order_0}) and (\ref{order_large}). From (\ref{order_0}) we learn that $r
\leq t$; and since $s \geq 0$ this is compatible with the assertion in
(\ref{PDE_reg}). In order to investigate the asymptotic behavior of
(\ref{order_large}) if $|\al| \geq 1$ we consider the cases $s = 0$ and $s > 0$
separately.
\begin{description}
\item{$s = 0$:} We recall that $r = \min(r,1) - \sig < 1$ and we have to show
that (\ref{order_large}) is $O(\eps^{r - |\al|})$. Combination of
(\ref{prod_bound}) (note that $|\al| \geq l \geq 1$) with the three possible
cases $O(1)$, $O(\log(1/\eps)$, $O(\eps^{t-l})$ of the growth rate of
$|b_\eps^{(l)}|$ directly yields an upper bound of the form
$O(\eps^{\min(t,1)-|\al|} (\log(1/\eps))^k)$ (where $k \leq l_0 + m_0 +1$).
Since the logarithmic factor is dominated by $\eps^{-\sig}$, for any $\sig >
0$, the assertion is proved.
\item{$s > 0$:} Now $r = \min(t,s+1)$ can be any real number and we have to go
through all cases relating the possible values of $|\al|$ and $r$.
\begin{trivlist}
\item[$|\al| < r$:] Since $1 \leq l \leq |\al|$ the factor
$|b_\eps^{(l)}|$ is $O(1)$, and in (\ref{der_struc2}), (\ref{prod_bound}) we
find $|\la| + |\mu| = |\al| - l < s$, which in turn yields $l_1 = m_1 = l_0 =
m_0 = 0$. Therefore we have an overall bound $O(1)$. \item[$|\al| = r$:] If $l
= |\al|$ then $|\la| + |\mu| = 0$ and hence $l_0 + m_0 = 0$ in
(\ref{prod_bound}) which means $O(1)$ for this part. The factor
$|b_\eps^{(l)}|$ gives at most $O(\log(1/\eps))$. If $l < |\al|$ then $l < t$
and so $|b_\eps^{(l)}|$ is $O(1)$. Since $|\la| + |\mu| \leq s$ we deduce $l_0
+ m_0 \leq 1$ and (\ref{prod_bound}) ensures an overall logarithmic bound.
Hence we have an upper bound of logarithmic order in both (sub)subcases.
\item[$|\al| > r$:] We note that $s > 0$ implies $(\log(1/\eps))^{l_0 + m_0}
\eps^{(l_0 + m_0)s} = O(1)$ whatever the value of $l_0 + m_0 \geq 0$.
Therefore the first case in (\ref{prod_bound}), $l_1 + m_1 \geq 1$, always
yields a bound $O(\eps^{s + l - |\al|})$.
\begin{trivlist}
\item[$l < t$:] The
$b$-dependent factor in (\ref{order_large}) is $O(1)$ and both cases in
(\ref{prod_bound}) are dominated by $O(\eps^{s + 1 - |\al|})$.
\item[$l = t$:]
$|b_\eps^{(t)}|$ contributes a logarithmic factor. The first case in
(\ref{prod_bound}) then gives $O(\eps^{s + t - |\al|})$ of which the part
$\eps^s$ can be used to suppress this logarithmic factor; hence a bound is
$O(\eps^{t - |\al|})$. The second case in (\ref{prod_bound}) yields an overall
bound which is some power of $\log(1/\eps)$ and therefore clearly dominated by
$\eps^{r - |\al|}$.
\item[$l > t$:] Here $|b_\eps^{(l)}| = O(\eps^{t - l})$.
Adding the factor according to the first case in (\ref{prod_bound}) then gives
a bound $O(\eps^{t - |\al| + s}) = O(\eps^{t - |\al|})$. On the other hand,
using the second line in (\ref{prod_bound}) provides an overall bound
$O(\eps^{t - l}\cdot (\log(1/\eps))^{l_0 + m_0})$. If $l < |\al|$ splitting
off $\eps^{t - |\al|}$ leaves an additional positive $\eps$-power to
compensate for the logarithmic term. If $l = |\al|$ we can again reason, like
in earlier cases, that $l_0 + m_0 = 0$. So, all branches of this
(subsub)subcase lead to a bound $O(\eps^{t - |\al|})$.
\end{trivlist}
Collecting the results of all (sub)subcases we have established the asymptotic
upper bound $O(\eps^{\min(t,s+1) - |\al|})$ of (\ref{order_large}).
\end{trivlist}
\end{description}
\end{proof}

The previous theorem indicates that we may expect a seismic wave to be about
one degree smoother than the irregular medium variation if the source is
prepared appropriately.
In principle this would enable one to deduce from measurements of the wave an
upper bound of the (global) medium regularity: first, estimate a strict upper
bound of the wave's Zygmund regularity $r$ via wavelet analysis of the data;
then the medium regularity cannot be better than $r-1$.

\paragraph{Acknowledgement} I thank Maarten de Hoop for initiating this line of
research and having gone through the very first steps in joint work with me
(\cite{HdH:01c}). He explained to me the relevance in applications and it was
his idea that combination of wavelet analysis with regularization in the
Colombeau setup could lead to new insights. However, the current paper would
not have been written without the continuous encouragement, valuable criticism
(in the best sense), and many suggestions for improvements by Michael
Oberguggenberger. This work was done while employed in his project P14576-MAT
by the Austrian Science Fund (FWF). 

\section*{Appendix: Characterization of Zygmund regularity via
continuous wavelet transform}

The proof to be presented below is a destillation of methods and basic
setups drawing from a variety of sources. We briefly sketch the basics of
these as a preparation.

\begin{trivlist}

\item[(i)]
Zygmund classes can alternatively be defined by a discrete Littlewood-Paley
decomposition (cf.\ \cite{Meyer:98}), also called dyadic resolution (e.g.,
in \cite{Triebel:II}). Let $\vphi_0 = \vphi$ and for $j\in\N$ put $
\vphi_j(\xi) = \int_{2^{j-1}}^{2^j} \psi(\xi/t)\, dt/t =
     \vphi(2^{-j}\xi) - \vphi(2^{-j+1}\xi)$.
 We have $\vphi_{j+1}(\xi) = \vphi_j(\xi/2)$ and the support of $\vphi_j$ ($j
 \geq 1$) is contained in the annulus $2^{j-1} \leq |\xi| \leq 2^{j+1}$. By
 construction, the family $(\vphi_j)_{j\geq 0}$ is a dyadic partition of unity:
 $\sum_{j=0}^\infty \vphi_j(\xi) = 1$.
 Similarly, the equation $\sum_j \vphi_j(D) u = u$ holds with convergence of
 the  series in $\S'$.

\item[(ii)]
The classical H\"older-Zygmund spaces can also be considered as the special
cases $B^s_{\infty,\infty}(\R^n)$  in Triebel's family of
Besov-Hardy-Sobolev-type spaces (cf.~\cite[Chapter 2, in particular
2.6.5/(1)]{Triebel:II}). These spaces are defined, for any $s\in\R$, by
$ B^s_{\infty,\infty} := \{u\in\S' \mid  \norm{u}{B^s_{\infty,\infty}} :=
        \sup_{j\geq 0} 2^{j s} \linf{\vphi_j(D)u} < \infty \}$.
The definition is independent of the particular choice of $\vphi$  (cf.\
\cite[2.3.2]{Triebel:II}).

\item[(iii)]  Both families of spaces, $B^s_{\infty,\infty}$ as well as
$\Zin^s$, are
realizations of the classical H\"older-Zygmund spaces when $s > 0$.
Therefore we clearly have $B^s_{\infty,\infty} = \Zin^s$ in this case. In
fact, equality holds for all real $s$: By \cite[2.3.8]{Triebel:I} (resp.\
\cite[Prop.8.6.6]{Hoermander:97}), for any $r\in\R$ the operators
$(1-\Delta)^{r/2}$ (resp.\ $(1-\Delta)^{-r/2}$) on $\S'$ map
$B^s_{\infty,\infty}$ (resp.\ $\Zin^{s-r}$) isomorphically into
$B^{s-r}_{\infty,\infty}$ (resp.\ $\Zin^s$); therefore we obtain
$  B^s_{\infty,\infty} = \Zin^s \quad \forall s\in\R$
with equivalent norms $\norm{u}{B^s_{\infty,\infty}}$ and $\inorm{u}{s}$.
We refer to these spaces as \emph{Zygmund spaces of regularity $s$}. In
particular, we deduce that the definition of $\Zin^s$ is independent of the
choice of $\vphi$.

\item[(iv)] In Meyer's book (cf.\ \cite[Chapter 3]{Meyer:98}) the
  H\"older-Zygmund spaces are
treated as special cases of Bony's two-microlocal spaces $C^{s,s'}_{x_0}$
(where $s > 0$, $s' = 0$, $x_0$ arbitrary). In fact, it is this point of
view which is underlying the proof of the characterization via the
('continuous') wavelet transform given in the following.

\end{trivlist}

\paragraph{Proof of Theorem \ref{wl_char}}
Recall that $\S_0(\R^n)$ is the subspace of $\S(\R^n)$ consisting of
functions with vanishing moments of all orders. Throughout the proof we
will make use of the following fact which will allow us to balance
vanishing moment conditions with regularity properties in occurring
convolutions.
\begin{lemma}\label{mom_reg} If $f\in\S$ with moments up to order $m-1$ vanishing
then one can find functions $f_\al\in\S$ ($|\al| = m$) such that
\[
    f = \sum_{|\al| = m} \d^\al f_\al.
\]
If, in addition, $f\in\S_0$ the functions $f_\al$ can be chosen to be in $\S_0$.
\end{lemma}
This can be shown by adapting the proof of \cite[Section 2.6, Lemma
12]{Meyer:92}. 

Concerning the notation of various constants in the estimates to follow we will
use the generic letter $C$, with subscripts if we want to indicate dependence on
certain parameters.

\emph{Part (i):} Applying the above lemma to $g$ we have $g_\eps = \eps^m
\sum_{|\al| = m} \d^\al ((g_\al)_\eps)$ and obtain
\[
    W_g u (.,\eps) = u * \ovl{\check{g}_\eps} =
        (-1)^m \eps^m \sum_{|\al| = m} (\d^\al u) * (\ovl{\check{g}_\al})_\eps .
\]
Since $\d^\al u \in \Zin^{s-m}$ with $s-m < 0$ and $g_\al\in\S$ we have reduced
the proof of (\ref{Wg_estimate}) to the task of estimating $\linf{u * g_\eps}$
where $u\in\Zin^s$ with $s<0$ and $g\in\S$.

If $T > 1$ arbitrary then $u = \vphi(D/T)u + \int_T^\infty \psi(D/T)u \, dt/t$
(with $\S'$-convergence) and therefore we have for any $\eps > 0$ fixed
 \begin{equation} \label{conv_dec}
    u * g_\eps = \vphi(D/T)u * g_\eps + \int_T^\infty \psi(D/t)u * g_\eps \, dt/t .
 \end{equation}
Let $T \geq 1/\eps \geq T/2$ and estimate the two terms in (\ref{conv_dec})
separately.

Since $\vphi(D/T)u = \vphi(D)u + \int_1^T \psi(D/T)u \, dt/t$ we deduce
(recalling that we may assume $s < 0$)
\begin{multline*}
    \linf{\vphi(D/T)u} \leq \linf{\vphi(D)u} + \int_1^T \linf{\psi(D/t)u} \, dt/t \\
    \leq C \Big( 1 + \int_1^T t^{-s} \, dt/t \Big) \leq C T^{-s} \leq 2^{-s} C \eps^s.
\end{multline*}
Therefore we obtain
 \begin{equation}
    \linf{\vphi(D/T)u * g_\eps} \leq \linf{\vphi(D/T)u} \lone{g_\eps} \leq C \eps^s.
 \end{equation}

To estimate the integrand in the second term of (\ref{conv_dec}) we assume $t
\geq T$ and choose
\newcommand{\tpsi}{\tilde{\psi}} $\tpsi\in\D$ with $\tpsi = 0$ near $0$ and $\tpsi = 1$ on
$\supp(\psi)$. It follows that $\F^{-1}\tpsi\in\S_0$ and $\psi(D/t)u * g_\eps =
\psi(D/t)u * \tpsi(D/t) g_\eps$.

Choose $r\in\N$ such that $r+s > 0$ and apply Lemma \ref{mom_reg} to obtain
functions $\tpsi_\al$, $|\al|=r$, satisfying $\F^{-1}\tpsi_\al\in\S_0$ and
$\F^{-1}\tpsi = \sum_{|\al|=r} D^\al \F^{-1}\tpsi_\al = \sum_{|\al|=r}
\F^{-1}(\xi^\al \tpsi_\al)$. Then $\tpsi(D/t) g_\eps
= t^{-r}\eps^{-r}  \sum_{|\al|=r} \tpsi_\al(D/t) (D^\al g)_\eps$ and since
$\linf{\psi(D/t)u} \leq C t^{-s}$ we have the estimate
\begin{multline*}
    \linf{\psi(D/t)u * g_\eps} \leq
    t^{-r} \eps^{-r} \sum_{|\al| =r} \linf{\psi(D/t)u * \tpsi_\al(D/t) (D^\al g)_\eps} \\
    \leq t^{-r} \eps^{-r} \sum_{|\al| =r} \linf{\psi(D/t)u}
        \lone{\tpsi_\al(D/t) (D^\al g)_\eps} \\
        \leq C t^{-(r+s)} \eps^{-r}
            \max\limits_{|\al|=r} \lone{\tpsi_\al(D/t) (D^\al g)_\eps}.
\end{multline*}
We show that the appearing $\L^1$-norms have bounds independent of $t$ and
$\eps$.

Writing $\tpsi_\al(D/t) (D^\al g)_\eps$ explicitly as a convolution and
rescaling by $t$ via substitution of the integration variable we have
\[
    \tpsi_\al(D/t) (D^\al g)_\eps(x) = \eps^{-n} \int \F^{-1}(\tpsi_\al)(y)
        (D^\al g)(\frac{x}{\eps}-\frac{y}{t\eps}) \, dy.
\]
For any $l$, the second factor in the integrand is bounded by $C_l
(1+|\frac{x}{\eps}-\frac{y}{t\eps}|^2)^{-l} \leq 2^l C_l (1 +
|\frac{x}{\eps}|^2)^{-l} (1 + |\frac{y}{t\eps}|^2)^l \leq 2^l C_l (1 +
|\frac{x}{\eps}|^2)^{-l} (1 + |y|^2)^l$ since $t\eps \geq T \eps \geq 1$.
Assuming $l > n/2$ and  integrating also over $x$ we finally obtain a bound for
$\lone{\tpsi_\al(D/t) (D^\al g)_\eps}$ of the form $C_{l,\al} \int \eps^{-n} (1
+ |\frac{x}{\eps}|^2)^{-l}\, dx = C_{l,\al} \int (1 + |z|^2)^{-l}\, dz$ which is
indeed independent of $t$ and $\eps$.

Taking the maximum of all bounds over $|\al| = r$ we arrive at the conclusion
that for all $t \geq T$
\[
    \linf{\psi(D/t)u * g_\eps} \leq C \eps^{-r} t^{-(r+s)}.
\]
If $R > T$ arbitrary then
 \begin{multline*}
    \linf{\int_T^R \psi(D/t)u * g_\eps \, dt/t}
        \leq \int_T^R \linf{\psi(D/t)u * g_\eps} \, dt/t \\
    \leq C \eps^{-r} \int_T^R t^{-(r+s)-1}\, dt =
    \frac{C\eps^{-r}}{r+s} (T^{-(r+s)} - R^{-(r+s)}).
 \end{multline*}
When $R\to \infty$ this upper bound tends to $\frac{C}{r+s} (\eps T)^{-r}\,
T^{-s} \leq C \eps^s$. This completes the proof of (\ref{Wg_estimate}).

\emph{Part (ii):}

\begin{lemma}\label{series_lemma}
Let $r>0$ and $k\in\N$ such that $k > r$. Assume that $h_j$ ($j\in\N_0$) is a
sequence of functions in $\Con^k(\R^n)$ with the property that there is $B > 0$
such that for all $\be\in\N_0^n$ with $|\be| \leq k$
 \begin{equation} \label{der_bounds}
    \linf{\d^\be h_j} \leq B\, 2^{j |\be|}.
 \end{equation}
Then the infinite series
 \begin{equation}
    h(x) := \sum_{j=0}^\infty 2^{-jr} h_j(x)
 \end{equation}
converges uniformly and defines an element in $\Zin^r(\R^n)$.
\end{lemma}
\begin{proof} Since $\linf{h_j} \leq B$ for all $j$ the series is
absolutely and uniformly convergent and defines a continuous bounded function
$h(x)$. Hence it is immediate that $\linf{\vphi(D)h} \leq \lone{\F^{-1}\vphi}
\linf{h}$. It remains to estimate $\linf{t^r \psi(D/t)h}$ for all $t \geq 1$. We
start by picking $q\in\N_0$ such that $2^q \leq t < 2^{q+1}$ and split the
necessary summation according to
\begin{multline*}
    | t^r \psi(D/t)h(x)| \leq \sum_{j=0}^\infty 2^{-jr} t^r |\psi(D/t)h_j(x)| \\
    = \sum_{j=0}^{q-1} 2^{-jr} t^r |\psi(D/t)h_j(x)| +
        \sum_{j=q}^\infty 2^{-jr} t^r |\psi(D/t)h_j(x)| =: S_1(x) + S_2(x).
\end{multline*}
The terms in $S_2$ can be estimated as follows
\begin{multline*}
    2^{-jr} t^r |\psi(D/t)h_j(x)| \leq
         (\frac{t}{2^q})^r 2^{-r(j-q)} \lone{\F^{-1}\psi} \linf{h_j} \\
    \leq 2^r C' B 2^{-r(j-q)} = C 2^{-r(j-q)}
\end{multline*}
and hence $S_2(x)$ is dominated uniformly by a convergent geometric series.

To find a bound for $S_1(x)$ we apply Lemma \ref{mom_reg} and rewrite
$\psi(D/t)$, as with $\tpsi$ in the proof of part (i), in the form $\psi(D/t) =
t^{-k} \sum_{|\al|=k} \psi_\al(D/t) D^\al$. Hence
\[
    |\psi(D/t)h_j(x)| \leq t^{-k} \sum_{|\al|=k} \lone{\F^{-1}\psi_\al} \linf{D^\al h_j}
    \leq t^{-k} C_\psi B 2^{jk} = C' t^{-k} 2^{jk}
\]
and we obtain
\begin{multline*}
    S_1(x) \leq C' \sum_{j=0}^{q-1} t^{r-k} 2^{j(k-r)} \leq
        C' 2^{-q(k-r)} \sum_{j=0}^{q-1} (2^{(k-r)})^j\\
    = C' 2^{-q(k-r)} \frac{2^{q(k-r)}-1}{2^{k-r}-1} \leq C.
\end{multline*}
Since $t\geq 1$ was arbitrary and the constants in the estimates are independent
of $q$ the lemma is proved.
\end{proof}

\begin{lemma}\label{Zyg_lemma} If $W\in\OM(\R^n\times\R_+)$ and satisfies (\ref{Wg_estimate}),
with $W$ substituted for $W_g u$, then
\[
    \int_0^1 W(.,\eps) * g_\eps \frac{d\eps}{\eps} \in  \Zin^s(\R^n).
\]
\end{lemma}
\begin{proof}
We show that the limit of $u^{(N)} := \int_{2^{-N}}^1 W(.,\eps)*g_\eps\, d\eps
/\eps$, as $N\to\infty$, defines an element in $\Zin^s$. As used already in part
(i) Lemma \ref{mom_reg} implies $g_\eps = \eps^m \sum_{|\al|=m}
\d^\al((g_\al)_\eps)$ and hence
\[
    u^{(N)} = \sum_{|\al|=m} \d_x^\al \Big( \underbrace{ \int_{2^{-N}}^1 \eps^{m-1}
        W(.,\eps)*(g_\al)_\eps}_{u^{(N)}_\al} \Big) \, d\eps.
\]
For any $\al\in\N_0^n$ with $|\al|=m$ the map $\d^\al\colon  \Zin^{s+m} \to
\Zin^s$ is continuous, hence it suffices to prove convergence of $u^{(N)}_\al$
in $\Zin^{s+m}$ (as $N\to\infty$) for each such $\al$.

By a dyadic subdivision of the interval $[2^{-N},1]$ we find a corresponding
series representation of $u^{(N)}_\al$ in the form
\begin{multline*}
    u^{(N)}_\al(x) =
    \sum_{j=0}^{N-1} \int_{2^{-j-1}}^{2^{-j}} \eps^{m-1}
        W(.,\eps) * (g_\al)_\eps (x)\, d\eps\\
    = \sum_{j=0}^{N-1} 2^{-j(m+s)} \int_{1/2}^1
    2^{js} W(.,2^{-j}\eta) * (g_\al)_{2^{-j}\eta} (x)\, \eta^{m-1} \, d\eta \\
    =: \sum_{j=0}^{N-1} 2^{-j(m+s)} v_{j,\al}(x)
\end{multline*}
where we have changed the variable $\eps = 2^{-j} \eta$. Note that $\linf{2^{js}
W(.,2^{-j}\eta)}$ $\leq C \eta^s \leq C$  independent of $j$. Therefore the
sequence $v_{j,\al}$ satisfies the condition (\ref{der_bounds}) of Lemma
\ref{series_lemma} for any $k\in\N$ with $k > m+s > 0$  since
\begin{multline*}
    |\d^\ga v_{j,\al}(x)| \leq
    \int_{1/2}^1 |2^{js} W(.,2^{-j}\eta) * (\d^\ga g_\al)_{2^{-j}\eta} (x)|\,
    2^{j|\ga|}\, \eta^{m - |\ga| - 1} \, d\eta \\
    \leq 2^{j|\ga|} \linf{2^{js} W(.,2^{-j}\eta)} \lone{\d^\ga g_\al}
    \int_{1/2}^1 \eta^{m - |\ga| - 1} \, d\eta = C_{\ga,\al}  2^{j|\ga|}.
\end{multline*}
Application of Lemma \ref{series_lemma} completes the proof.
\end{proof}

\begin{lemma} \label{smooth_lemma} Let $W\in\OM(\R^n\times\R_+)$ then
\[
    \int_1^\infty W(.,\eps) * g_\eps \frac{d\eps}{\eps} \in  \Cinf(\R^n).
\]
\end{lemma}
\begin{proof} Let $R > 1$ and put $v_R = \int_1^R W(.,\eps)*g_\eps\, d\eps/\eps$. Then $v_R$
is smooth, temperate, and converges weakly to some $v\in\S'$ as $R\to\infty$
(cf.\ (\ref{synthesis})).


Clearly, any derivative $\d^\al v_R$ converges to $\d^\al v$ then. But letting
the derivative fall on the factor $g_\eps$ inside the integral defining $v_R$
produces additional factors $\eps^{-|\al|}$. When $|\al|$ is large enough to
compensate for the polynomial growth of $W(y,\eps)$ with respect to $\eps$ this
ensures absolute convergence of the classical integral. Hence for all $|\al|$
sufficiently large $\d^\al v$ is smooth, yielding that $v$ itself is smooth.
\end{proof}

To finish the proof of part (ii) we apply (\ref{synthesis}) together with
(\ref{p_inversion}) and obtain, with some polynomial $p$,
\[
    u = \int_0^1 W_g u (.,\eps) * g_\eps \frac{d\eps}{\eps}
        + \int_1^\infty W_g u (.,\eps) * g_\eps \frac{d\eps}{\eps} + p.
\]
The second term is smooth by Lemma \ref{smooth_lemma} and the first term is of
Zygmund regularity $s$ by Lemma \ref{Zyg_lemma}. It follows that $u$ differs
from an element in $\Zin^s$ only by some smooth function.

This completes the proof of Theorem \ref{wl_char}.

\bibliographystyle{abbrv}
\bibliography{gueMO}

\end{document}